\documentclass[12pt]{amsart}
\usepackage{amsfonts,amssymb,amsmath,setspace,color}

\newtheorem{theorem}{Theorem}[section]
\newtheorem{lemma}[theorem]{Lemma}
\newtheorem{proposition}[theorem]{Proposition}
\newtheorem{corollary}[theorem]{Corollary}
\theoremstyle{definition}
\newtheorem{definition}[theorem]{Definition}

\newtheorem{remark}[theorem]{Remark}

\theoremstyle{claim}

\begin{document}

\title[Maps preserving common zeros]{\bf Maps preserving common zeros between subspaces
of vector-valued continuous functions}

\author{Luis Dubarbie}

\address{Departamento de Matem\'aticas, Estad\'istica y
Computaci\'on, Facultad de Ciencias, Universidad de Cantabria,
Avenida de los Castros s/n, E-39071, Santander, Spain.}

\email{luis.dubarbie@gmail.com}

\thanks{Research supported by the Spanish Ministry
of Science and Education (MTM2006-14786) and by a predoctoral grant from
the University of Cantabria and the Goverment of Cantabria.}



\renewcommand{\theequation}{\Alph{equation}}

\begin{abstract}
For metric spaces $X$ and $Y$, normed spaces $E$ and $F$, and certain subspaces
$A(X,E)$ and $A(Y,F)$ of vector-valued continuous functions,
we obtain a complete characterization of linear and bijective maps
$T:A(X,E)\rightarrow A(Y,F)$ preserving common zeros, that is,
maps satisfying the property
\setcounter{equation}{15}
\begin{equation}\label{dub}
Z(f)\cap Z(g)\neq \emptyset \Longleftrightarrow Z(Tf)\cap Z(Tg)\neq \emptyset
\end{equation}
for any $f,g\in A(X,E)$, where $Z(f)=\{x\in X:f(x)=0\}$. Moreover, we provide some
examples of subspaces for which the automatic continuity of linear bijections
having the property (\ref{dub}) is derived.
\end{abstract}

\maketitle

\section{Introduction}

In the recent paper \cite{LT}, the authors deal with linear and bijective maps
$T:A(X,E)\rightarrow A(Y,F)$ such that
\addtocounter{equation}{9}
\begin{equation}\label{leu}
\cap_{i=1}^{k}Z(f_{i})\neq \emptyset \Longleftrightarrow \cap_{i=1}^{k}Z(Tf_{i})\neq \emptyset
\end{equation}
for any $k\in \mathbb{N}$ and any sequence $(f_{i})_{i=1}^{k}$ in $A(X,E)$.
They obtain a general description of vector space
isomorphisms $T:A(X,E)\rightarrow A(Y,F)$ satisfying (\ref{leu}), where $A(X,E)$ and $A(Y,F)$ denote
certain subspaces of $C(X,E)$ and $C(Y,F)$, respectively, for completely regular spaces $X$ and $Y$,
and Hausdorff topological vector spaces $E$ and $F$. More precisely, when $X$ and $Y$ are realcompact
spaces, they prove that every vector space isomorphism $T:C(X,E)\rightarrow C(Y,F)$ having the property (\ref{leu})
can be written as a weighted composition operator and show that $X$ and $Y$ are homeomorphic
and $E$ and $F$ are isomorphic. They also consider such maps defined between spaces of vector-valued
differentiable functions and obtain an analogous result.

\smallskip

\renewcommand{\theequation}{\arabic{equation}}
\setcounter{equation}{0}

A related problem was studied in \cite{CRX}, where the authors consider,
for compact Hausdorff spaces $X$ and $Y$, and a Banach lattice $E$,
Riesz isomorphisms $T:C(X,E)\rightarrow C(Y,\mathbb{R})$ that preserve
non-vanishing functions. They showed that $X$ and $Y$ are homeomorphic
and $E$ and $\mathbb{R}$ are Riesz isomorphic. Later, this result has been generalized
for different kinds of Banach lattices (see \cite{MCX,EO,EO1}). Finally,
in \cite{CCW}, it was proved that $X$ and $Y$ are homeomorphic and
$E$ and $F$ are Riesz isomorphic whenever $T:C(X,E)\rightarrow C(Y,F)$
is a Riesz isomorphism preserving non-vaninshing functions, $X$ and $Y$
are compact Hausdorff spaces, and $E$ and $F$ are arbitrary Banach lattices.
Recently, the Lipschitz version of this topic has been obtained in \cite{JMV},
providing a complete characterization of vector lattice isomorphisms
$T:A(X,E)\rightarrow A(Y,F)$ preserving non-vanishing functions, where $X$ and $Y$ are
compact metric spaces, $E$ and $F$ are Banach lattices, and $A(X,E)$ and $A(Y,F)$ are vector
sublattices of $\mathrm{Lip}(X,E)$ and $\mathrm{Lip}(Y,F)$, respectively, that separate
and join points uniformly.

\smallskip

In this paper, given $X$ and $Y$ metric spaces, $E$ and $F$ normed spaces,
and $A(X,E)$ and $A(Y,F)$ certain subspaces of $C(X,E)$ and $C(Y,F)$, respectively,
we deal with linear bijections $T:A(X,E)\rightarrow A(Y,F)$ having the property (\ref{dub}). Concretely,
we obtain a complete description of such maps and prove the existence of a homeomorphism
between $X$ and $Y$. Finally, we give some examples of subspaces of vector-valued continuous
functions for which the automatic continuity of maps preserving common zeros is derived.

\section{Preliminaries and definitions}

Let $X$ be a metric space and $E$ be a $\mathbb{K}$-normed space, where $\mathbb{K}$ stands for
the field of real or complex numbers. $C(X,E)$ will denote the set of all
$E$-valued continuous functions defined on $X$. If $E=\mathbb{K}$,
we set $C(X):=C(X,E)$.

\smallskip

Given a function $f:X\rightarrow E$, we define the \emph{zero set}
associated to $f$ as $Z(f):=\{x\in X: f(x)=0\}$ and the \emph{cozero set} of $f$ as
$\mathrm{coz}(f):=X\setminus Z(f)$. The constant scalar-valued function defined on $X$ and
taking the value $1$ will be denoted by ${\bf 1}_{X}$. Given $f\in C(X)$ and $g\in C(X,E)$,
$f\cdot g$ stands for their pointwise product.
Finally, if $A$ is a subset of $X$, then $int(A)$ and $cl(A)$ denote the interior and the
closure of $A$ in $X$, respectively.

\smallskip

From now on, $X$ and $Y$ will be metric spaces
and $E$ and $F$ will be $\mathbb{K}$-normed spaces. Besides, we will consider
vector subspaces $A(X,E)$ of $C(X,E)$ (respectively $A(Y,F)\subset C(Y,F)$)
and subrings $A(X)$ of $C(X)$ (respectively $A(Y)\subset C(Y)$) satisfying
the following properties:
\begin{enumerate}
\item $A(X,E)$ contains non-vanishing functions.
\item $A(X,E)$ is an $A(X)$-module, that is, given $f\in A(X)$ and $g\in A(X,E)$, then $f\cdot g\in A(X,E)$.
\item For each $x_{0}\in X$, there exists $f_{x_{0}}\in A(X,E)$ such that
$$Z(f_{x_{0}})=\{x_{0}\}.$$
\item $A(X)$ is completely regular, that is, given a point $x_{0}$ in $X$ and
a closed subset $C$ of $X$ such that $x_{0}\notin C$, there exists $\varphi\in A(X)$
satisfying $\varphi(x_{0})=1$ and $\varphi\equiv 0$ on $C$.
\end{enumerate}

\begin{definition}
A linear map $T:A(X,E)\rightarrow A(Y,F)$ is said to \emph{preserve common zeros} if
$$
Z(f)\cap Z(g)\neq \emptyset \Longleftrightarrow Z(Tf)\cap Z(Tg)\neq \emptyset
$$
for any $f,g\in A(X,E)$.
\end{definition}

\begin{remark}\label{prop}
It is easily seen that, given a linear and bijective map $T:A(X,E)\rightarrow A(Y,F)$
preserving common zeros, then
$$
Z(f)\neq \emptyset \Longleftrightarrow Z(Tf)\neq \emptyset
$$
for any $f\in A(X,E)$.
\end{remark}

As it has been pointed out in the previous section, Remark \ref{prop} shows that maps
preserving common zeros between subspaces of vector-valued continuous functions are
strongly related with maps preserving non-vanishing functions.

\begin{definition}
A linear map $T:A(X,E)\rightarrow A(Y,F)$ is said to be \emph{separating} if
$\mathrm{coz}(Tf)\cap \mathrm{coz}(Tg)=\emptyset$
whenever $f,g\in A(X,E)$ satisfy $\mathrm{coz}(f)\cap \mathrm{coz}(g)=\emptyset.$
Moreover, we will say that $T$ is \emph{biseparating} if it is bijective and both $T$ and $T^{-1}$
are separating.
\end{definition}

\section{Maps preserving common zeros}

In this section, we provide a complete description of
bijective linear maps preserving common zeros between subspaces of vector-valued
continuous functions that satisfy the properties given above.
Concretely, the main result we obtain is the following.

\begin{theorem}\label{imp}
Let $T:A(X,E)\rightarrow A(Y,F)$ be a linear and bijective map
that preserves common zeros. Then there exist a homeomorphism $h:Y\rightarrow X$,
$\mathbb{K}$-normed subspaces $E_{y}\subset E$ and $F_{y}\subset F$, and a map
$Jy:E_{y}\rightarrow F_{y}$ linear and bijective for each $y\in Y$, such that
$$
Tf(y)=Jy(f(h(y)))
$$
for all $f\in A(X,E)$ and $y\in Y$.
\end{theorem}

The rest of the section is devoted to prove the previous result.
From now on, to simplify the notation, we assume that $T$ is a bijective linear map
preserving common zeros from $A(X,E)$ into $A(Y,F)$. Besides,
$\Phi_{X}$ and $\Phi_{Y}$ will denote non-vanishing functions in $A(X,E)$ and $A(Y,F)$, respectively.

\begin{proposition}\label{hola}
Let $f\in A(X,E)$ be such that $Z(f)=\{x_{0}\}$. Then there exists a unique
$y_{0}\in Y$ with $y_{0}\in Z(Tf)$.
\end{proposition}

\begin{proof}
Obviously $Z(Tf)\neq \emptyset$ by Remark \ref{prop}.\\
Suppose now that there exist distinct points $y_{1},y_{2}\in Y$ such that $y_{1},y_{2}\in Z(Tf)$.
As we are assuming that $A(Y)$ is completely regular, there exists
$\varphi \in A(Y)$ such that $\varphi(y_{1})=1$ and $\varphi(y_{2})=0$.
Next, taking $\Phi_{Y}\in A(Y,F)$, we define the functions
$
g:=\varphi\cdot \Phi_{Y}
$
and
$
l:=({\bf 1}_{Y}-\varphi)\cdot \Phi_{Y},
$
which belong to $A(Y,F)$ since it is an $A(Y)$-module. Immediately,
$Z(T^{-1}g)\neq \emptyset$ and $Z(T^{-1}l)\neq \emptyset$ by Remark \ref{prop}.
Now, it is clear that $y_{2}$ belongs to $Z(g)\cap Z(Tf)$,
and then $Z(T^{-1}g)\cap Z(f)\neq \emptyset$. Also, since $Z(f)=\{x_{0}\}$, we deduce that
$x_{0}\in Z(T^{-1}g)$. Similarly, $Z(T^{-1}l)\cap Z(f)\neq \emptyset$ because
$y_{1}\in Z(l)\cap Z(Tf)$, and then $x_{0}\in Z(T^{-1}l)$. Therefore,
$$
0=T^{-1}g(x_{0})+T^{-1}l(x_{0})=T^{-1}\Phi_{Y}(x_{0}),
$$
and consequently $Z(\Phi_{Y})\neq \emptyset$, which is absurd.\\
Finally, assume that there exists $f_{1}\in A(X,E)$ not equal to $f$ such that
$Z(f_{1})=\{x_{0}\}$. Then, by above reasoning, $Z(Tf_{1})$ contains a unique point,
and, since $Z(Tf)\cap Z(Tf_{1})\neq \emptyset$, we can conclude that $Z(Tf_{1})=Z(Tf)=\{y_{0}\}$.
\end{proof}

\begin{remark}\label{46}
The previous proposition lets us define a map $k:X\rightarrow Y$ sending each point
$x\in X$ such that $Z(f)=\{x\}$ for some $f\in A(X,E)$ to the point $y\in Y$
satisfying $Z(Tf)=\{y\}$.
\end{remark}

\begin{lemma}\label{36}
Let $f\in A(X,E)$ and $x_{0}\in X$ be such that $f(x_{0})=0$. Then $Tf(k(x_{0}))=0$.
\end{lemma}

\begin{proof}
Take $f_{x_{0}}\in A(X,E)$ such that $Z(f_{x_{0}})=\{x_{0}\}$.
Then, clearly $Z(Tf_{x_{0}})\cap Z(Tf)\neq \emptyset$.
Now, by Proposition \ref{hola} and Remark \ref{46}, $Z(Tf_{x_{0}})=\{k(x_{0})\}$,
and consequently $k(x_{0})$ belongs to $Z(Tf)$.
\end{proof}

\begin{remark}
Since $T$ is a bijective map that preserves common zeros, we can consider the map
$T^{-1}:A(Y,F)\rightarrow A(X,E)$ which also preserves common zeros. As a consequence,
applying Proposition \ref{hola} and Remark \ref{46}, we obtain a map
$h:Y\rightarrow X$ associated to $T^{-1}$ for which we have the following result.
\end{remark}

\begin{lemma}\label{09}
Let $g\in A(Y,F)$ and $y_{0}\in Y$ be such that $g(y_{0})=0$. Then $T^{-1}g(h(y_{0}))=0$.
\end{lemma}

\begin{theorem}\label{homeo}
$h:Y\rightarrow X$ is a homeomorphism.
\end{theorem}

\begin{proof}
Firstly, we will see that $h$ is a bijective map. It is enough to prove that $k\equiv h^{-1}$,
and therefore, we check that $(k\circ h)(y)=y$ for all $y\in Y$ and $(h\circ k)(x)=x$ for all $x\in X$.
Suppose, on the contrary, that there exists $y_{0}\in Y$ such that $k(h(y_{0}))\neq y_{0}$. Now,
since $A(Y)$ is completely regular, we can consider $\varphi \in A(Y)$ satisfying $\varphi(k(h(y_{0})))=1$
and $\varphi(y_{0})=0$, and define the function $g:=\varphi\cdot \Phi_{Y}\in A(Y,F)$.
It is clear that $g(y_{0})=0$ and, by Lemma \ref{09}, $T^{-1}g(h(y_{0}))=0$. Applying Lemma \ref{36}, we conclude that
$[T(T^{-1}g)](k(h(y_{0})))=0$, that is, $g(k(h(y_{0})))=0$, which is impossible by definition of $g$.
The equality $(h\circ k)(x)=x$ for all $x\in X$ can be seen in a similar way.\\
Secondly, we will prove that $h$ and $k$ are continuous functions. Suppose that $h$ is not continuous.
Then there exists a sequence $(y_{n})$ in $Y$ converging to $y_{0}\in Y$ such that
$(h(y_{n}))$ does not converge to $h(y_{0})$.
Therefore, there exists an open neighborhood $U$ of $h(y_{0})$ such that $h(y_{n})\notin U$
for infinitely many $n\in \mathbb{N}$. Since we are assuming that $A(X)$ is completely regular,
we can take $\varphi \in A(X)$ such that $\varphi(h(y_{0}))=1$ and
$\varphi\equiv 0$ on $X\setminus U$. Finally, we define $f:=\varphi\cdot \Phi_{X}\in A(X,E)$.
For infinitely many $n\in \mathbb{N}$, $f(h(y_{n}))=0$. Next, taking into account Lemma \ref{36},
we obtain that $Tf(k(h(y_{n})))=0$ for infinitely many $n\in \mathbb{N}$, or equivalently, $Tf(y_{n})=0$.
As a consequence, $Tf(y_{0})=0$, and applying Lemma \ref{09}, $[T^{-1}(Tf)](h(y_{0}))=0$, that is, $f(h(y_{0}))=0$,
in contradiction with the definition of the function $f$.
The fact that $k$ is continuous is proved similarly.
\end{proof}

\begin{definition}\label{810}
For each $y\in Y$, we consider the $\mathbb{K}$-normed subspaces
$$
E_{y}:=\{f(h(y)):f\in A(X,E)\}
$$
and
$$
F_{y}:=\{g(y):g\in A(Y,F)\}
$$
of $E$ and $F$, respectively.
\end{definition}

\begin{remark}
Since $A(X,E)$ and $A(Y,F)$ contain non-vanishing functions, for every $y\in Y$,
$E_{y}$ and $F_{y}$ are non-trivial $\mathbb{K}$-normed subspaces of $E$ and $F$, respectively.
\end{remark}

\begin{definition}
For each $y\in Y$, we define the map $Jy:E_{y}\rightarrow F_{y}$ as $Jy(\textmd{e}):=Tf(y)$
for all $\textmd{e}\in E_{y}$, where $f\in A(X,E)$ satisfies $f(h(y))=\textmd{e}$.
\end{definition}

\begin{remark}
It is not difficult to check that $Jy$ is well defined for each $y\in Y$.
Fix $y_{0}\in Y$ and consider $\textmd{e}\in E_{y_{0}}$. If there exist
distinct functions $f_{1}$ and $f_{2}$ in $A(X,E)$ with
$f_{1}(h(y_{0}))=\textmd{e}=f_{2}(h(y_{0}))$, applying Lemma \ref{36}
we deduce that $T(f_{1}-f_{2})(y_{0})=0$, and consequently $Tf_{1}(y_{0})=Tf_{2}(y_{0})$,
as we wanted to see.
\end{remark}

\begin{lemma}\label{J}
$Jy$ is linear and bijective for all $y\in Y$.
\end{lemma}

\begin{proof}
Fix $y_{0}\in Y$. It is clear that $Jy_{0}$ is linear. Next,
we see that $Jy_{0}$ is injective. Let $\textmd{e}\in E_{y_{0}}$ be such that
$
Jy_{0}(\textmd{e})=0.
$
By definition, $Tf(y_{0})=0$ for some $f\in A(X,E)$ such that $f(h(y_{0}))=\textmd{e}$.
On the other hand, by Lemma \ref{09}, $[T^{-1}(Tf)](h(y_{0}))=0$, which implies that $\textmd{e}=0$.
Finally, we prove that $Jy_{0}$ is an onto map.
Fix $\textmd{f}\in F_{y_{0}}$. Therefore, there exists $g\in A(Y,F)$ such that $g(y_{0})=\textmd{f}$.
Besides, since $T$ is onto, $g=Tf$ for some $f\in A(X,E)$. Then, defining $\textmd{e}:=f(h(y_{0}))$,
it is obvious that $\textmd{e}\in E_{y_{0}}$ and
$$
Jy_{0}(\textmd{e})=Tf(y_{0})=g(y_{0})=\textmd{f}.
$$
Since this can be done for every $y\in Y$, the lemma is proved.
\end{proof}

\begin{proof}[Proof of Theorem \ref{imp}]
Assume that $T$ is a linear and bijective map that preserves common zeros.
Then, for each $f\in A(X,E)$ and $y\in Y$, clearly
$T$ has the desired representation by definition of $Jy$. Moreover, the existence of
a homeomorphism $h$ between $Y$ and $X$ was proved in Theorem \ref{homeo}, and finally,
in Lemma \ref{J} has been shown that each $Jy$ is linear and bijective.
\end{proof}

\begin{remark}
A simple example shows that the homeomorphism $h$ between $Y$ and $X$
obtained in Theorem \ref{imp} is somehow determined by the subspaces $F_{y}$ of $F$.
Let $X=\{a,b,c\}$, $Y=\{1,2,3\}$, $E=F=\mathbb{K}^{4}$, $A(X)=C(X)$ and $A(Y)=C(Y)$.
Now, consider three subspaces $E^{a}$, $E^{b}$, and $E^{c}$ of $E$ with $\mathrm{dim}(E^{a})=1$,
$\mathrm{dim}(E^{b})=2$, and $\mathrm{dim}(E^{c})=3$, and define
$$
A(X,E):=\left\{f\in C(X,E): f(a)\in E^{a}, f(b)\in E^{b}, f(c)\in E^{c} \right\}.
$$
Finally, for each $i\in \{1,2,3\}$, take $F^{i}\subset F$ such that $\mathrm{dim}(F^{i})=i$
and define
$$
A(Y,F):=\left\{g\in C(Y,F): g(i)\in F^{i}, i\in Y \right\}.
$$
Obviously $A(X,E)$, $A(Y,F)$, $A(X)$, and $A(Y)$ satisfy the properties given in Section 2.
Now, if we suppose that there exists a linear bijection preserving common zeros between
$A(X,E)$ and $A(Y,F)$, by Definition \ref{810},
$F_{i}=\{g(i): g\in A(Y,F)\}=F^{i}$ for $i=1,2,3$.
On the other hand, we have proved that, for each $i\in \{1,2,3\}$,
there exists a linear and bijective map $J_{i}:E_{i}\rightarrow F_{i}$
where $E_{i}=\{f(h(i)): f\in A(X,E)\}$, so $h(1)=a$, $h(2)=b$, and $h(3)=c$.

However, in general, the homeomorphism $h:Y\rightarrow X$ is not unique. Consider $X=Y=\mathbb{N}$,
$E=F=\mathbb{K}^{3}$, $A(X)=C(X)$, $A(Y)=C(Y)$, and define
$$
A(X,E):=\left\{f\in C(X,E): f(2n-1)\in E^{2n-1}, f(2n)\in E^{2n},  n\in \mathbb{N} \right\}
$$
where $E^{2n-1}$ and $E^{2n}$ are subspaces of $E$ with $\mathrm{dim}(E^{2n-1})=1$ and
$\mathrm{dim}(E^{2n})=2$, and
$$
A(Y,F):=\left\{g\in C(Y,F):\ g(2n-1)\in F^{2n-1},\ g(2n)\in F^{2n}, n\in \mathbb{N} \right\}
$$
for $F^{2n-1}, F^{2n} \subset F$ with $\mathrm{dim}(F^{2n-1})=2$ and $\mathrm{dim}(F^{2n})=1$.
Following the previous reasoning, it is clear that $h(2n-1)$ is an even number
and $h(2n)$ is an odd number for each $n\in \mathbb{N}$.
\end{remark}

\begin{corollary}\label{bi}
Let $T:A(X,E)\rightarrow A(Y,F)$ be a linear and bijective map preserving common zeros.
Then $T$ is biseparating.
\end{corollary}

\begin{proof}
To prove that $T$ is biseparating, we will see that both $T$ and its inverse are separating.
Let $f,g\in A(X,E)$ be such that $\mathrm{coz}(Tf)\cap \mathrm{coz}(Tg)\neq \emptyset$. Then
there exists $y_{0}\in Y$ satisfying $Tf(y_{0})\neq 0$ and $Tg(y_{0})\neq 0$.
By representation of $T$ given in Theorem \ref{imp}, it is obvious that $Jy_{0}(f(h(y_{0})))\neq 0$ and
$Jy_{0}(g(h(y_{0})))\neq 0$. Therefore, since $Jy_{0}$ is injective,
$h(y_{0})\in \mathrm{coz}(f)\cap \mathrm{coz}(g)$, which implies
that $T$ is a separating map. Analogously it can be proved that $T^{-1}$ is also
separating.
\end{proof}

\section{Applications I}

This section is devoted to prove that every linear bijection that preserves common zeros
between spaces of vector-valued differentiable functions defined on open subsets of $\mathbb{R}^{p}$,
$p\in \mathbb{N}$, is automatically continuous. Firstly, we recall some basic facts about
such spaces.

\smallskip

Let $X$ be an open subset of $\mathbb{R}^{p}$, $p\in \mathbb{N}$, and let $E$ be a
$\mathbb{K}$-Banach space. Then $C^{n}(X,E)$ consists of the space of $E$-valued functions
defined on $X$ whose partial derivaties up to order $n$ exist and are continuous.
Besides, when $X$ is also bounded, $C^{n}(\overline{X},E)$ denotes
the subspace of $C^{n}(X,E)$ of functions whose partial derivaties up to order $n$ can
be continuously extended to the boundary of $X$.

\smallskip

From now on, we will assume that we are in one of the following contexts.

\smallskip

{\bf Context 1.} $X$ and $Y$ will be open subsets of $\mathbb{R}^{p}$ and $\mathbb{R}^{q}$,
respectively ($p,q\in \mathbb{N}$), and $E$ and $F$ will be $\mathbb{K}$-Banach spaces.
$A(X,E)=C^{n}(X,E)$, $A(Y,F)=C^{m}(Y,F)$, $A(X)=C^{n}(X)$, and $A(Y)=C^{m}(Y)$ for any $n,m\in \mathbb{N}$.

\smallskip

{\bf Context 2.} $X$ and $Y$ will be bounded open subsets of $\mathbb{R}^{p}$ and $\mathbb{R}^{q}$,
respectively ($p,q\in \mathbb{N}$), with the property that $int(cl(X))=X$ and $int(cl(Y))=Y$. $E$ and $F$
will be $\mathbb{K}$-Banach spaces. $A(X,E)=C^{n}(\overline{X},E)$, $A(Y,F)=C^{m}(\overline{Y},F)$,
$A(X)=C^{n}(\overline{X})$, and $A(Y)=C^{m}(\overline{Y})$ for any $n,m\in \mathbb{N}$.

\smallskip

This means that when we refer to $A(X,E)$, $A(Y,F)$, $A(X)$, and $A(Y)$,
we suppose that all of them are included at the same time in one of the
previous contexts.

\smallskip

Since our purpose is to prove the automatic continuity of maps preserving common zeros
between spaces of differentiable functions, we will endow such spaces with
some natural topology.

\begin{definition}
A locally convex topology in $A(X,E)$ is said to be \emph{compatible with the pointwise convergence}
if the following conditions are satisfied:
\begin{enumerate}
\item $A(X,E)$ endowed with such topology is a Banach space.
\item If $(f_{n})$ is a sequence in $A(X,E)$ converging to $0$, then $(f_{n}(x))$ converges to $0$
for all $x\in X$.
\end{enumerate}
\end{definition}

It is not difficult to see that $A(X,E)$ and $A(X)$ satisfy the properties
given in Section 2 (analogously for $A(Y,F)$ and $A(Y)$). Therefore, Theorem \ref{imp}
and Corollary \ref{bi} can be applied in the previous contexts. Consequently,
taking into account \cite[Theorem 5.1]{A}, it is immediate to deduce the following result.

\begin{theorem}
Let $T:A(X,E)\rightarrow A(Y,F)$ be a linear and bijective map preserving common zeros.
If we assume that $A(X,E)$ and $A(Y,F)$ are endowed with any topologies which are compatible
with the pointwise convergence, then $T$ is continuous.
\end{theorem}

\section{Applications II}

Finally, we provide more examples of subspaces of vector-valued continuous
functions to which the previous results can be applied. Besides, we prove that,
under some conditions, every bijective linear map preserving common zeros
between such subspaces is continuous.

\smallskip

As in the previous section, we will assume that we are in one of the two following contexts.

\smallskip

{\bf Context 1.} $X$ and $Y$ are bounded complete metric spaces,
$E$ and $F$ are $\mathbb{K}$-Banach spaces, $A(X,E)=\mathrm{Lip}(X,E)$
denotes the space of all bounded $E$-valued Lipschitz functions on $X$,
and $A(X)=\mathrm{Lip}(X)$ (respectively $A(Y,F)=\mathrm{Lip}(Y,F)$ and
$A(Y)=\mathrm{Lip}(Y)$) (see \cite{W}).

\smallskip

{\bf Context 2.} $X$ and $Y$ are compact subsets of the real line,
$E$ and $F$ are $\mathbb{K}$-Banach spaces, $A(X,E)=AC(X,E)$ stands for the space
of absolutely continuous functions defined on $X$ and taking values in $E$,
and $A(X)=AC(X)$ (respectively $A(Y,F)=AC(Y,F)$ and $A(Y)=AC(Y)$)
(see \cite[Section 18]{HS}).

\smallskip

For the two previous contexts, it is easy to check that  $A(X,E)$, $A(Y,F)$,
$A(X)$, and $A(Y)$ satisfy the properties given in Section 2. Finally,
the continuity of linear and bijective maps that preserve common zeros
between such subspaces can be deduced.

\begin{theorem}
Let $T:A(X,E)\rightarrow A(Y,F)$ be a bijective linear map preserving common zeros.
Then $T$ is continuous if one of the following conditions is satisfied:
\begin{enumerate}
\item $E$ and $F$ has the same finite dimension.
\item $Y$ has no isolated points.
\end{enumerate}
\end{theorem}

\begin{proof}
Taking into account Corollary \ref{bi}, it is clear that $T$ is a biseparating map.
Therefore, if we consider Context 1, the result is obtained applying \cite[Corollaries 5.11 and 5.12]{AD}.
For the second context, by Theorems 3.6 and 4.8 in \cite{DU}, we conclude that $T$ is continuous.
\end{proof}

\proof[Acknowledgements] The author wishes to thank the suggestions given
by Professor J. Araujo.

\end{document}